\documentclass[12pt,psamsfonts]{amsart}

\newtheorem{anyprop}{Anyprop}[section]

\newtheorem{theorem}[anyprop]{Theorem}
\newtheorem{lemma}[anyprop]{Lemma}
\newtheorem{proposition}[anyprop]{Proposition}

\theoremstyle{definition}

\newtheorem{definition}[anyprop]{Definition}

\newtheorem{remark}[anyprop]{Remark}

\newtheorem{theoremintro}{Theorem}

\newtheorem{definitionintro}{Definition}

\newcommand{\NN}{\mathbb{N}}

\newcommand{\PP}{\mathbb{P}}

\newcommand  {\shF}     {\mathcal{F}}

\newcommand  {\shK}     {\mathcal{K}}

\newcommand  {\shL}     {\mathcal{L}}

\newcommand  {\shS}     {\mathcal{S}}
\newcommand  {\shT}     {\mathcal{T}}

\newcommand  {\shQ}     {\mathcal{Q}}

\newcommand  {\fom}     {\mathfrak{m}}

\newcommand  {\dual}    {\vee}

\newcommand  {\Ext}     {\operatorname{Ext}}

\newcommand  {\lra}     {\longrightarrow}

\renewcommand{\O}       {\mathcal{O}}

\newcommand  {\Proj}    {\operatorname{Proj}}

\newcommand  {\ra}      {\rightarrow}

\newcommand  {\rk}    {\operatorname{rk}}

\newcommand  {\Spec}    {\operatorname{Spec}}

\newcommand  {\Syz}     {\operatorname{Syz}}

\newcommand{\komdots}{ , \ldots , }
\newcommand{\plusdots}{ + \ldots + }

\newcommand{\leqdots}{ \leq \ldots \leq }

\newcommand{\subsetdots}{ \subset \ldots \subset }

\newcommand{\numiii}{\renewcommand{\labelenumi}{(\roman{enumi})}}

\newcommand{\length} {\lambda}

\theoremstyle{remark}

\numberwithin{equation}{section}

\usepackage{amscd}
\usepackage{amssymb}

\def\mydate{\number\day\space\ifcase\month \or January\or February\or March\or April\or May\or
June\or July\or August\or September\or October\or November\or
December\fi \space\number\year}

\newcommand{\mumu}{\mu_{HK}}

\begin{document}

\title[Hilbert-Kunz criterion]
{A Hilbert-Kunz criterion  for solid closure in dimension two (characteristic zero)}

\author[Holger Brenner]{Holger Brenner}
\address{Mathematische Fakult\"at, Ruhr-Universit\"at Bochum, 
               44780 Bochum, Germany}
\email{Holger.Brenner@ruhr-uni-bochum.de}





\begin{abstract}
Let $I$ denote a homogeneous $R_+$-primary ideal in a two-dimensional normal standard-graded domain over an algebraically closed field of characteristic zero.
We show that a homogeneous element $f$ belongs to the solid closure $I^*$ if and only
if $e_{HK}(I) = e_{HK}((I,f))$, where $e_{HK}$ denotes the (characteristic zero) Hilbert-Kunz multiplicity
of an ideal. This provides a version in characteristic zero of the well-known Hilbert-Kunz criterion
for tight closure in positive characteristic.
\end{abstract}

\maketitle


\section*{Introduction}

Let $(R, \fom)$ denote a local Noetherian ring or an $\NN$-graded algebra
of dimension $d$ of positive characteristic $p$.
Let $I$ denote an $\fom$-primary ideal, and set $I^{[q]}= (f^q: f\in I)$ for a prime power $q=p^{e}$.
Then the Hilbert-Kunz function of $I$ is given
by
$$ e \longmapsto \length ( R/I^{[p^{e}]}) \, ,$$
where $\length $ denotes the length.
The Hilbert-Kunz multiplicity of $I$ is defined as the limit
$$ e_{HK} (I)= \lim_{e \ra \infty}  \, \length (R/I^{[p^{e}]}) /p^{ed} \, .$$
This limit exists as a positive real number, as shown by Monsky in \cite{monskyhilbertkunz}.
It is an open question whether this number is always rational.

The Hilbert-Kunz multiplicity is related to the theory of tight closure.
Recall that the tight closure of an ideal $I$ in a Noetherian ring
of characteristic $p$ is by definition the ideal 
$$I^* \!\!=\!\! \{f \in R: \exists c \mbox{ not in any minimal prime} : cf^q \in I^{[q]}
\mbox{ for almost all } q=p^{e}\} \, .$$
For an analytically unramified and formally equidimensional local ring $R$
the equation 
$e_{HK}(I) = e_{HK} (J) $ holds if and only if $I^* =J^*$ holds true for ideals $I \subseteq J$
(see \cite[Theorem 5.4]{hunekeapplication}).
Hence $f \in I^*$ if and only if $e_{HK}(I) = e_{HK} ( (I,f))$.
This is the Hilbert-Kunz criterion for tight closure in positive characteristic.

The aim of this paper is to give a characteristic zero version of this relationship between
Hilbert-Kunz multiplicity and tight closure for $R_+$-primary
homogeneous ideals in a normal two-dimensional graded domain $R$.
There are several notions for tight closure in characteristic zero,
defined either by reduction to positive characteristic or directly. We will work with the notion of solid closure (see \cite{hochstersolid}). In dimension two, the containment in the solid closure
$f \in (f_1 \komdots f_n)^*$ means that the open subset $D( \fom) \subset \Spec A$
is not an affine scheme,
where $A=R[T_1 \komdots T_n]/(f_1T_1 \plusdots f_nT_n+f)$ is the so-called forcing algebra,
see \cite[Proposition 1.3]{brennertightproj}.

The definition of the Hilbert-Kunz multiplicity in positive characteristic does not suggest at first sight
an analogous notion in characteristic zero.
However, a bridge is provided by
the following result of \cite{brennerhilbertkunz}, which gives an explicit formula for the Hilbert-Kunz
multiplicity and proves its rationality in dimension two
(the rationality of the Hilbert-Kunz multiplicity for the maximal ideal was also obtained independently in \cite{trivedihilbertkunz}).

\begin{theoremintro}
Let $R$ denote a two-dimensional standard-graded normal domain over an algebraically closed field of positive
characteristic, $Y= \Proj R$.
Let $I=(f_1 \komdots f_n)$ denote a homogeneous $R_+$-primary ideal
generated by homogeneous elements $f_i$ of degree $d_i, i=1 \komdots n$.
Then the Hilbert-Kunz multiplicity of the ideal $I$  equals
$$e_{HK}(I)=\frac{ \deg (Y)}{2}( \sum_{k=1}^t r_k \nu_k^2 - \sum_{i=1}^n d_i^2) \, .$$ 
\end{theoremintro}

Here the numbers $r_k$ and $\nu_k$ come from the strong Harder-Narasimhan filtration of
the syzygy bundle $\Syz(f_1^q \komdots f_n^q)(0)$
given by the short exact sequence
$$0 \lra \Syz(f_1^q \komdots f_n^q)(0) \lra \bigoplus_{i=1}^n \O(-q d_i) \stackrel{f_1^q \komdots f_n^q}{\lra}
\O_Y \lra 0 \, .$$
This syzygy bundle is a locally free sheaf on the smooth projective curve $Y= \Proj R$,
and its strong Harder-Narasimhan filtration
is a filtration $\shS_1 \subsetdots \shS_t = \Syz(f_1^q \komdots f_n^q)(0)$
such that the quotients $\shS_k/\shS_{k-1}$ are strongly semistable, meaning that every Frobenius pull-back is semistable. Such a filtration exists for $q$ big enough by a theorem of Langer,
\cite[Theorem 2.7]{langersemistable}.
Then we set $r_k = \rk (\shS_k/\shS_{k-1})$ and $\nu_k=- \mu( \shS_k/\shS_{k-1}) /q \deg (Y)$,
where $\mu$ denotes the slope.

To define the Hilbert-Kunz multiplicity in characteristic zero we now take the right hand side of the above formula as our model.

\begin{definitionintro}
Let $R$ denote a two-dimensional normal standard-graded $K$-domain over an algebraically closed field $K$
of characteristic zero. Let $I=(f_1 \komdots f_n)$ 
be a homogeneous $R_+$-primary ideal given by homogeneous ideal generators $f_i$ of degree $d_i$.
Let $\shS_1 \subsetdots \shS_t = \Syz(f_1 \komdots f_n)(0)$ denote the Harder-Narasimhan filtration of
the syzygy bundle on $Y= \Proj R$,
set $\mu_k =\mu( \shS_k/\shS_{k-1}) $ and $r_k = \rk ( \shS_k/\shS_{k-1}) $.
Then the Hilbert-Kunz multiplicity of $I$ is by definition
$$e_{HK}(I)= \frac{ \deg (Y)}{2}( \sum_{k=1}^t r_k (\frac{\mu_k}{\deg(Y)})^2 - \sum_{i=1}^n d_i^2) 
= \frac{ \sum_{k=1}^t r_k \mu_k^2 - \deg(Y)^2 \sum_{i=1}^n d_i^2}{2  \deg (Y)} \, .$$ 
\end{definitionintro}

It is easy to show that this definition does not depend
on the chosen ideal generators and is therefore an invariant of the ideal,
see \cite[Proposition 4.9]{brennerhilbertkunz}.
With this invariant we can in fact give the following Hilbert-Kunz criterion
for solid closure in characteristic zero in dimension two (see Theorem \ref{hilbertkunzsolid}):

\begin{theoremintro}
\label{hilbertkunzsolidintro}
Let $K$ denote an algebraically closed field of characteristic zero,
let $R$ denote a standard-graded two-dimensional normal $K$-domain.
Let $I$ be a homogeneous $R_+$-primary ideal and let $f$ denote a homogeneous element.
Then $f$ is contained in the solid closure,
$f \in I^*$, if and only if $e_{HK} (I)= e_{HK}((I,f))$.
\end{theoremintro}

To prove this theorem it is convenient to consider more generally for a locally free sheaf $\shS$
on a smooth projective curve $Y$ the expression
$$\mumu(\shS) = \sum_{k=1}^t  r_k \mu_k^2 \, ,$$
where $r_k$ and $\mu_k$ are the ranks and the slopes of the semistable quotient
sheaves in the Harder-Narasimhan filtration of $\shS$. We call this number the Hilbert-Kunz slope
of $\shS$. With this notion the Hilbert-Kunz multiplicity of an ideal $I=(f_1 \komdots f_n)$
is related to the Hilbert-Kunz slope
of the syzygy bundle by
$$e_{HK}((f_1 \komdots f_n))
=  \frac{1}{2  \deg (Y)} \big( \mumu(\Syz(f_1 \komdots f_n)(0)) - \mumu (\bigoplus_{i=1}^n \O(-d_i)) \big) \, .$$
With this notion we will in fact prove the following theorem, which implies
Theorem \ref{hilbertkunzsolidintro} (see Theorem \ref{hilbertkunzaffinecriterion}).

\begin{theoremintro}
\label{hilbertkunzaffinecriterionintro}
Let $Y$ denote a smooth projective curve over an algebraically closed field of characteristic $0$.
Let $\shS$ denote a locally free sheaf on $Y$ and let $c \in H^1(Y, \shS)$
denote a cohomology class given rise to the extension $0 \ra \shS \ra \shS' \ra \O_Y \ra 0$
and the affine-linear torsor $\PP(\shS'^\dual) - \PP(\shS^\dual)$.
Then $\PP(\shS'^\dual) - \PP(\shS^\dual)$ is an affine scheme if and only if
$\mumu(\shS') < \mumu (\shS)$.
\end{theoremintro}

\section{The Hilbert-Kunz slope of a vector bundle}

We recall briefly some notions for locally free sheaves (or vector bundles), see \cite{huybrechtslehn}
or \cite{hardernarasimhan}.
Let $Y$ denote a smooth projective curve over an algebraically closed field
and let $\shS$ denote a locally free sheaf of rank $r$. Then
$\deg (\shS)= \deg (\bigwedge^r \shS)$ is called the degree of $\shS$
and $\mu(\shS)= \deg(\shS)/r$ is called the slope of $\shS$.
If $\mu(\shT) \leq \mu(\shS)$ holds for every locally free subsheaf $\shT \subseteq \shS$,
then $\shS$ is called semistable.
In general there exists the so-called Harder-Narasimhan filtration.
This is a filtration of locally free subsheaves
$\shS_1 \subsetdots \shS_t = \shS$ such that
the quotient sheaves $\shS_k/\shS_{k-1}$ are semistable locally free sheaves
with decreasing slopes $\mu_1 > \ldots > \mu_t$.
The Harder-Narasimhan filtration is uniquely determined by these properties.
$\shS_1$ is called the maximal destabilizing subsheaf, $\mu_1 =\mu_{\rm max}(\shS)$
is called the maximal slope of $\shS$
and $\mu_t= \mu_{\rm min}(\shS)$ is called the minimal slope of $\shS$.
If $\shS \ra \shT$ is a non-trivial sheaf homomorphism,
then $\mu_{\rm min}(\shS) \leq \mu_{\rm max}(\shT)$.

We begin with the definition of the Hilbert-Kunz slope of $\shS$.

\begin{definition}
Let $\shS$ denote a locally free sheaf on a smooth projective curve
over an algebraically closed field of characteristic $0$.
Let $\shS_1 \subsetdots \shS_t=\shS$ denote the Harder-Narasimhan filtration of $\shS$,
set $r_k = \rk (\shS_k/\shS_{k-1})$ and $\mu_k = \mu(\shS_k/\shS_{k-1})$.
We define the Hilbert-Kunz slope of $\shS$ by
$$\mumu (\shS) = \sum_{k=1}^t r_k \mu_k^2 = \sum_{k=1}^t \frac{ \deg( \shS_k/\shS_{k-1}  )^2}{ r_k}\, .$$
\end{definition}

The only justification for considering this number is Theorem \ref{hilbertkunzsolid} below.
We gather together some properties of this notion in the following proposition.

\begin{proposition}
\label{hilbertkunzslopeproperties}
Let $\shS$ denote a locally free sheaf on a smooth projective curve over an algebraically closed field
of characteristic $0$. Then the following hold true.

\numiii

\begin{enumerate}

\item
If $\shS$ is semistable, then $\mumu(\shS) =\deg(\shS)^2/\rk(\shS)$.

\item
Let $\shT \subset \shS$ denote a locally free subsheaf occurring in the Harder-Narasimhan filtration
of $\shS$.
Then $\mumu(\shS) =\mumu(\shT) + \mumu(\shS/\shT)$.

\item
We have $\mumu(\shS \oplus \shT) = \mumu(\shS) + \mumu(\shT)$.

\item
$\mumu(\shS)= \mumu(\shS^\dual)$.

\item
Let $\shL$ denote an invertible sheaf. Then
$$\mumu (\shS \otimes \shL) = \mumu(\shS) + 2 \deg(\shS) \deg (\shL) + \rk(\shS) \deg(\shL)^2 \, .$$

\item
Let $\varphi: Y' \ra Y$ denote a finite morphism between smooth projective curves of degree $n$.
Then $\mumu( \varphi^*(\shS)) = n^2 \mumu(\shS)$.
\end{enumerate}

\end{proposition}
\proof
(i) and (ii) are clear from the definition.
(iii). The maximal destabilizing subsheaf of $\shS \oplus \shT$
is either $\shS_1 \oplus 0$, $0 \oplus \shT_1$ or $\shS_1 \oplus \shT_1$.
Hence the result follows from (ii) by induction on the rank of $\shS \oplus \shT$.

(iv). Let $0= \shS_0 \subset \shS_1 \subsetdots \shS_t = \shS$ denote the Harder-Narasimhan filtration
of $\shS$.
Set $\shQ_k = \shS/\shS_k$.
This gives a filtration
$0 \subset \shQ_{t-1}^\dual \subsetdots \shQ_1^\dual \subset \shQ_0^\dual = \shS^\dual$.
From $0 \ra \shS_k /\shS_{k-1} \ra \shS/\shS_{k-1} \ra \shS/\shS_k \ra 0$
we get
$0 \ra \shQ_k^\dual \ra \shQ_{k-1}^\dual \ra \shQ^\dual_{k-1}/\shQ^\dual_k
\cong ( \shS_k/\shS_{k-1})^\dual
\ra 0$.
Hence the filtration is the Harder-Narasimhan filtration of $\shS^\dual$ and
the result follows from $\mu( \shQ^\dual_{k-1}/\shQ^\dual_k) = - \mu(\shS_k/\shS_{k-1})$.

(v).
The Harder-Narasimhan filtration of $\shS \otimes \shL$
is $\shS_1 \otimes \shL \subsetdots \shS_t \otimes \shL $ and
$\mu(\shS_k \otimes \shL/ \shS_{k-1} \otimes \shL )= \mu( (\shS_k/ \shS_{k-1}) \otimes \shL )
=\mu (\shS_k/ \shS_{k-1}) + \mu(\shL)$.
Therefore
\begin{eqnarray*}
\mumu(\shS \otimes \shL)
&=& \sum_{k=1}^t r_k \mu_k( \shS \otimes \shL)^2 \cr
&=& \sum_{k=1}^t r_k (\mu_k + \deg (\shL))^2 \cr
&=& \sum_{k=1}^t r_k (\mu_k^2 + 2 \mu_k \deg (\shL) +\deg(\shL)^2) \cr
&=& \mumu (\shS) + 2\deg(\shL) \sum_{k=1}^t r_k\mu_k  +\deg(\shL^2)\sum_{k=1}^t r_k \, .
\end{eqnarray*}
This is the stated result, since $\deg(\shS) = \sum_{k=1}^t r_k \mu_k$
and $\rk(\shS)= \sum_{k=1}^t r_k$.

(vi). The pull-back of a semistable sheaf under a separable morphism is again semistable,
and the pull-back of the Harder-Narasimhan filtration is the Harder-Narasimhan filtration
of $\varphi^* (\shS)$. Hence the result follows from $\deg(\varphi^*(\shS)) = n \deg (\shS)$.
\qed

\begin{lemma}
\label{semistableminimal}
The Hilbert-Kunz multiplicity of a locally free sheaf $\shS$ has the property that
$\mumu(\shS) \geq  \deg(\shS)^2/\rk(\shS)$, and equality holds if and only if $\shS$ is semistable.
\end{lemma}
\proof
We have to show that
$$\sum_{k=1}^t r_k \mu_k ^2 \geq \deg(\shS)^2/\rk(\shS)
= (r_1\mu_1 \plusdots r_t\mu_t)^2/(r_1 \plusdots  r_t) $$
or equivalently that
$$(r_1 \plusdots  r_t) (\sum_{k=1}^t r_k \mu_k ^2 ) \geq (r_1\mu_1 \plusdots r_t\mu_t)^2 \,.$$
The left hand side is
$ \sum_{k=1}^t r_k^2 \mu_k^2 + \sum_{i \neq k} r_i r_k \mu_k^2$ (we sum over ordered pairs),
and the right hand side is
$ \sum_{k=1}^t r_k^2 \mu_k^2 + \sum_{i \neq k} r_i r_k \mu_i \mu_k$.
Hence left hand minus right hand is
$$  \sum_{i \neq k} r_i r_k \mu_k^2 - \sum_{i \neq k} r_i r_k \mu_i \mu_k $$
So this follows from $0 \leq (\mu_i - \mu_k)^2 =\mu_i^2 + \mu_k^2 -2 \mu_i \mu_k$
for all pairs $i \neq k$. Equality holds if and only if $\mu_i= \mu_k$, but then $t=1$ and
$\shS$ is semistable.
\qed

\begin{remark}
Lemma \ref{semistableminimal} implies that the number
$\mumu(\shS) - \frac{ \deg(\shS)^2}{\rk(\shS)} \geq0$, and $=0$ holds exactly in the semistable case.
It follows from Proposition \ref{hilbertkunzslopeproperties} (v)
that this number is invariant under tensoring with an invertible sheaf.
\end{remark}

\begin{proposition}
Let $\shS$ and $\shT$ denote two locally free sheaves on $Y$.
Then
$$\mumu(\shS \otimes \shT)= \rk(\shT) \mumu(\shS) + \rk(\shS) \mumu(\shT) +2 \deg(\shS) \deg(\shT)    \, .$$
\end{proposition}
\proof
Let $r_i$, $\mu_i$, $i \in I$, and $r_j$, $\mu_j$, $j \in J$, ($I$ and $J$ disjoint)
denote the ranks and slopes occurring in the Harder-Narasimhan filtration of $\shS$ and $\shT$ respectively.
It is a non-trivial fact (in characteristic zero!) that the tensor product of two semistable bundle
is again semistable, see \cite[Theorem 3.1.4]{huybrechtslehn}.
From this it follows that the semistable quotients of the Harder-Narasimhan filtration of
$\shS \otimes \shT$ are given as $(\shS_i/ \shS_{i-1}) \otimes (\shT_j/\shT_{j-1})$
of rank $r_i \cdot r_j$ and slope $\mu_i+ \mu_j$.
Therefore the Hilbert-Kunz slope is
\begin{eqnarray*}
\mumu(\shS \otimes \shT) &=& \sum_{i,j} r_i r_j (\mu_i + \mu_j)^2 \cr
&=& \sum_{i,j} r_i r_j \mu_i^2 + \sum_{i,j}r_ir_j \mu_j^2 +2\sum_{i,j} r_ir_j \mu_i \mu_j \cr
&=& (\sum_{j} r_j )( \sum_i r_i \mu_i^2) + (\sum_{i} r_i )( \sum_j r_j \mu_j^2)
+2( \sum_i  r_i \mu_i)(\sum_j r_j \mu_j) \cr
&=& \rk(\shT) \mumu(\shS) + \rk(\shS) \mumu(\shT) +2 \deg(\shS) \deg(\shT) 
\end{eqnarray*}
\qed

\section{A Hilbert-Kunz criterion for affine torsors}
\label{}

In this section we consider a locally free sheaf $\shS$ on a smooth projective curve $Y$
together with a cohomology class $c \in H^1(Y, \shS) \cong \Ext(\O_Y, \shS)$. Such a class gives rise
to an extension
$0 \ra \shS \ra \shS' \ra \O_Y \ra 0$. Of course $\deg(\shS')= \deg(\shS)$ and
$\rk(\shS')= \rk (\shS) + 1$.
We shall investigate the relationship between $\mumu(\shS)$ and $\mumu(\shS')$.

\begin{lemma}
\label{extensionlemma}
Let $Y$ denote a smooth projective curve over an algebraically closed field.
Let $\shS$, $\shT$ and $\shQ$ denote locally free sheaves on $Y$.
Then the following hold.

\numiii

\begin{enumerate}

\item
Let $\varphi: \shT \ra \shS$ denote a sheaf homomorphism, $c \in H^1(Y, \shT)$ with corresponding
extension $\shT'$, let $\shS'$ denote the extension of $\shS$ corresponding to
$\varphi(c) \in H^1(Y, \shS)$.
Then there is a sheaf homomorphism $\varphi': \shT' \ra \shS'$ extending $\varphi$.

\item
Suppose that $0 \ra \shT \ra \shS \ra \shQ \ra 0$ is a short exact sequence,
and $c \in H^1(Y, \shT)$. Then $\shT' \subseteq \shS'$ and $\shS'/ \shT' \cong \shS/ \shT$.

\item
Suppose that $0 \ra \shT \ra \shS \ra \shQ \ra 0$ is a short exact sequence,
and $c \in H^1(Y, \shS)$. Then $\shS' \ra \shQ' \ra 0$ and $\shQ' \cong \shS'/\shT$.

\item
If $\shS$ is semistable of degree $0$ and $c \in H^1(Y,\shS)$,
then also $\shS'$ is semistable.
\end{enumerate}

\end{lemma}
\proof
The cohomology class $c$ is represented by the $\check{C}$ech cocycle
$\check{c} \in H^0(U_1 \cap U_2, \shS)$, where $Y= U_1 \cup U_2$ is an affine covering.
Then $\shS'$ arises from $\shS_1'= \shS|_{U_1} \oplus \O$
and $\shS_2'= \shS|_{U_2} \oplus \O$ by glueing
$\shS_1'|{U_1 \cap U_2} \cong \shS_2'|{U_1 \cap U_2}$ via
$(s,t) \mapsto (s+ t\check{c},t)$.
The natural mappings $\shT_i' \ra \shS_i'$, $i=1,2$, glue together to a morphism
$\shT' \ra \shS'$.
The injectivity and surjectivity transfer from $\varphi$ to $\varphi'$, since these are local properties.
(ii) and (iii) then follow from suitable diagrams.

(iv). Suppose that $\shF \subseteq \shS'$ is a semistable subsheaf of positive slope.
Then the induced mapping $\shF \ra \O$ is trivial and therefore $\shF \subseteq \shS$,
which contradicts the semistability of $\shS$.
\qed

\medskip
Let $\shS_1 \subsetdots \shS_t =\shS$ denote the Harder-Narasimhan filtration of $\shS$
and $c \in H^1(Y, \shS)$.
If the image of $c$ in $H^1(Y, \shS/\shS_{t-1})$ is zero, then $c$ stems
from a class $c_{t-1} \in H^1(Y, \shS_{t-1})$.
So we find inductively a class $c_n \in H^1(Y, \shS_n)$
mapping to $c$ and such that the image in $H^1(Y, \shS_n/\shS_{n-1})$ is not zero
(or $c$ itself is $0$).
This yields extensions $\shS_k'$ of $\shS_k$ for $k \geq n$.
It is crucial for the behavior of $\shS'$
whether $\mu(\shS_n/\shS_{n-1}) \geq 0$ or $< 0$.
The following Proposition deals with the case $\mu(\shS_n/\shS_{n-1}) \geq 0$.

\begin{proposition}
\label{positive}
Let
$\shS_1 \subsetdots \shS_t=\shS$
be the Harder-Narasimhan filtration of $\shS$
and let $c \in H^1(Y,\shS)$.
Let $n$ be such that the image of $c$ in $H^1(Y, \shS_k/\shS_{k-1})$ is $0$
for $k >n$ but such that the image in $H^1(Y,\shS_n/ \shS_{n-1})$ is $ \neq 0$.
Suppose that
$\mu (\shS_n/ \shS_{n-1})$ is $ \geq 0$.
Let $i$ be the biggest number such that $ \mu(\shS_i/\shS_{i-1}) \geq 0$
{\rm(}hence $n \leq i${\rm)}.

\numiii

\begin{enumerate}

\item
Suppose that $\mu_i > 0$. Then the Harder Narasimhan filtration of $\shS'$ is
$$\shS_1 \subsetdots  \shS_i \subset \shS_i' \subset \shS_{i+1}' \subsetdots \shS' \, .$$

\item
Suppose that
$\mu_i= 0$.
Then the Harder-Narasimhan filtration of $\shS'$ is
$$\shS_1 \subsetdots \shS_{i-1} \subset \shS_i' \subset \shS_{i+1}'
\subsetdots \shS' \, .$$
\end{enumerate}
\end{proposition}
\proof
(i).
The quotients of the filtration are $ \shS_k/ \shS_{k-1}$, $k \leq i$,
which have positive slope,
$\shS_i' /\shS_i \cong \O_Y$, and
$\shS_k' / \shS'_{k-1} \cong \shS_k /\shS_{k-1}$ (Lemma \ref{extensionlemma}(ii)) for $k >i$,
which have negative slope.
These quotients are all semistable and the slope numbers are decreasing.

(ii).
The quotients $\shS_k/\shS_{k-1}$ are semistable with decreasing positive slopes for $k=1 \komdots i-1$.
The quotients $\shS_k'/ \ \shS_{k-1}' \cong \shS_k /\shS_{k-1}$
are semistable with decreasing negative slopes for $k =i+1 \komdots  t$.
The quotient $\shS_i' /\shS_{i-1}$ is isomorphic to $(\shS_i/\shS_{i-1})'$
by Lemma \ref{extensionlemma}(iii), hence semistable of degree $0$ by
Lemma \ref{extensionlemma}(iv).
\qed

\medskip
In the rest of this section we study the remaining case, that $\mu( \shS_n/ \shS_{n-1}) <0$.
In this case it is not possible to describe the Harder-Narasimhan filtration
of $S'$ explicitly. However we shall see that in this case the Hilbert-Kunz slope of $\shS'$ is smaller
than the Hilbert-Kunz slope of $\shS$.
We need the following two lemmata.

\begin{lemma}
\label{slopecompare}
Let $\shT$ denote a locally free sheaf on $Y$ with Harder-Narasimhan filtration $\shT_k$,
$\mu_k = \mu(\shT_k/\shT_{k-1})$ and $r_k= \rk(\shT_k/\shT_{k-1})$.
Let
$$ (\tau_i) =( \mu_1 \komdots \mu_1, \mu_2 \komdots \mu_2, \mu_3 \komdots
\mu_{t-1}, \mu_t \komdots \mu_t)$$
denote the slopes where each $\mu_k$ occurs $r_k$-times.
Let $\shS \subseteq \shT$ denote a locally free subsheaf of rank $r$
and let $\sigma_i$, $i=1 \komdots r$ denote the corresponding numbers for $\shS$.
Then $\sigma_i \leq \tau_i$ for $i=1 \komdots r$.

Moreover, if $\shS$ is saturated {\rm(}meaning that the quotient sheaf is locally free{\rm)}
and if no subsheaf $\shS_j$ of the Harder-Narasimhan filtration of $\shS$ occurs
in the Harder-Narasimhan filtration of $\shT$,
then  $\sigma_i \leq \tau_{i+1}$ for $i= 1 \komdots r$.
\end{lemma}

\proof
Let $i$, $i=1 \komdots r$ be given and let $j$ be such that
$\rk (\shS_{j-1}) < i \leq \rk(\shS_j)$, hence $\sigma_i = \mu_j(\shS) = \mu(\shS_j/\shS_{j-1})$.
We may assume that $i= \rk(\shS_j)$.
Let  $k$ be such that $\rk(\shT_{k-1}) < i \leq \rk(\shT_k)$.
Therefore
$\shS_j \not \subseteq \shT_{k-1}$, and the induced morphism $\shS_j \ra \shT /\shT_{k-1}$ is not trivial.
Hence
$\sigma_i=\mu_j(\shS) =\mu_{\rm min}(\shS_j) \leq \mu_{\rm max} (\shT/\shT_{k-1}) = \mu_k(\shT)=\tau_i$.

Now suppose that $ \sigma_i > \tau_{i+1}$.
Then necessarily $\sigma_i > \sigma_{i+1}$ and $\tau_i > \tau_{i+1}$ by what we have already proven.
Therefore $i= \rk(\shS_j) = \rk(\shT_k)$.
If $\shS_j \subseteq \shT_k$, then they are equal, since both sheaves are saturated of the same rank,
but this is excluded by the assumptions.
Hence $\shS_j \not\subseteq \shT_k$ and $\shS_j \ra \shT /\shT_k$ is non-trivial.
Therefore $\sigma_i = \mu_{\rm min} (\shS_j) \leq \mu_{\rm max}(\shT/\shT_k)= \mu_{k+1}(\shT)= \tau_{i+1}$.
\qed

\begin{remark}
If the numbers $\tau_i$ are given as in the previous lemma, then
$\deg(\shT) = \sum_i \tau_i$ and $\mumu(\shT ) = \sum_{i} \tau_i^2$.
\end{remark}

\begin{lemma}
\label{numkrit}
Let $ \alpha_1 \leqdots \alpha_r$ and $\beta_1 \leqdots \beta_{r+1}$
denote  positive real numbers such that $\alpha_{i} \geq \beta_{i+1} $
for $i=1 \komdots r$ and $\sum_{i=1}^r \alpha_i = \sum_{i=1}^{r+1} \beta_i$.
Then $  \sum_{i=1}^{r+1} \beta_i^2 \leq \sum_{i=1}^r \alpha_i^2   $
and equality holds if and only if $\alpha_i = \beta_{i+1}$.
\end{lemma}
\proof
Let $\alpha_i = \beta_{i+1} +\delta_i$, $\delta_i \geq 0$.
From $\sum_{i=1}^r \alpha_i = \sum_{i=1}^r \delta_i + \sum_{i=1}^r \beta_{i+1}= \sum_{i=1}^{r+1} \beta_i$
we get $\beta_1 = \sum_{i=1}^r \delta_i$ ($\leq \beta_2 $).
The quadratic sums are
$$ \sum_{i=1}^{r} \alpha_i^2 = \sum_{i=2}^{r+1}    \beta_i^2 +  \sum_{i=1}^r \delta_i^2 +2 \sum_{i=1}^r \delta_i \beta_{i+1} $$
and
$$ \sum_{i=1}^{r+1} \beta_i^2 = ( \sum_{i=1}^r \delta_i)^2 + \sum_{i=2}^{r+1} \beta_i^2
= 2 \sum_{i< j }\delta_i \delta_j    +  \sum_{i=1}^r \delta_i^2 + \sum_{i=2}^{r+1} \beta_i^2 \, .$$
So we have to show that
$ \sum_{i <j } \delta _i \delta _j \leq \sum_{j=1}^r \delta_i \beta_{i+1}$.
But this is clear from
$\sum_{i < j} \delta_j \leq \sum_{j=1}^r \delta_j \leq \beta_2 \leq \beta_{i+1}$
for all $i= 1 \komdots r$.
Equality holds if and only if $\delta_i=0$.
\qed

\medskip
A cohomology class $H^1(Y, \shS)$ corresponds to a geometric $\shS$-torsor $T \ra Y$.
This is an affine-linear bundle on which $\shS$ acts transitively.
A geometric realization is given as $T = \PP(\shS'^\dual) - \PP(\shS^\dual)$.
The global cohomological properties of this torsor are related
to the Hilbert-Kunz slope in the following way.

\begin{theorem}
\label{hilbertkunzaffinecriterion}
Let $Y$ denote a smooth projective curve over an algebraically closed field of characteristic $0$.
Let $\shS$ denote a locally free sheaf on $Y$ and let $c \in H^1(Y, \shS)$
denote a cohomology class given rise to the extension $0 \ra \shS \ra \shS' \ra \O_Y \ra 0$
and the affine-linear torsor $\PP(\shS'^\dual) - \PP(\shS^\dual)$.
Then the following are equivalent.

\numiii

\begin{enumerate}

\item
There exists a locally free quotient $\varphi: \shS \ra \shQ \ra 0$
such that $\mu_{\rm max}(\shQ) <0$ and the image $\varphi(c) \in H^1(Y, \shQ)$ is non-trivial.

\item
The torsor $\PP(\shS'^\dual) - \PP(\shS^\dual)$ is an affine scheme.

\item
The Hilbert-Kunz slope drops, that is $\mumu(\shS') < \mumu(\shS)$.
\end{enumerate}

\end{theorem}
\proof
The equivalence (i) $\Leftrightarrow$ (ii) was shown in \cite[Theorem 2.3]{brennertightplus}.
The implication (iii) $\Rightarrow$ (i) follows from Proposition \ref{positive}:
for if (i) does not hold, then we are in the situation of
Proposition \ref{positive} that $\mu(\shS_n/\shS_{n-1}) \geq 0$.
The explicit description of the Harder-Narasimhan filtration
of $\shS'$ gives in both cases that $\mumu(\shS') = \mumu(\shS)$.

So suppose that (i) holds. This means that there exists a subsheaf
$\shS_n \subseteq \shS$ occurring in the Harder-Narasimhan filtration of $\shS$ such that $c$
stems from $c_{n} \in H^1(Y, \shS_{n})$ and such that
its image in $ H^1(Y, \shS/\shS_{n-1})$
is non-trivial with $ \mu_{\rm max} (\shS/\shS_{n-1}) =\mu(\shS_{n}/\shS_{n-1}) = \mu_n < 0$.

Let $\shT_1 \subsetdots \shT_t=\shS'$ denote the Harder-Narasimhan filtration
of $\shS'$ with slopes $\mu_k= \mu(\shT_k/\shT_{k-1})$ and ranks $r_k=\rk(\shT_k/\shT_{k-1})$.
Suppose that the maximal slope $\mu(\shT_1)$ is positive.
Then the induced mapping $\shT_1 \ra \shS'/\shS \cong \O_Y$ is trivial,
and $\shT_1 \subseteq \shS$.
This is then also the maximal destabilizing subsheaf of $\shS$, since 
$\mu_{\rm max} (\shS) \leq \mu_{\rm max}(\shS') =\mu(\shT_1)$.
Therefore $\mumu(\shS) = \mumu(\shT_1) + \mumu(\shS/\shT_1)$
and $\mumu(\shS') = \mumu(\shT_1) + \mumu(\shS'/\shT_1)$
by Proposition \ref{hilbertkunzslopeproperties}(ii).
Since $\shS'/\shT_1$ is the extension of $\shS/\shT_1$ defined by the image of the cohomology class in
$H^1(Y, \shS/\shT_1)$ (Lemma \ref{extensionlemma}(iii)),
we may mod out $\shT_1$. Note that this does not change the condition in (i).
Hence we may assume inductively that $\mu_{\rm max}(\shS) \leq 0$ and
$\mu_{\rm max}(\shS') \leq 0$.

Now suppose that $\shT_1$ has degree $0$. Again, if $\shT_1 \subseteq \shS$,
then this is also the maximal destabilizing subsheaf of $\shS$, and we can mod out $\shT_1$ as before.
So suppose that $\shT_1 \ra \O_Y$ is non-trivial. Then this mapping is surjective,
let $\shK \subset \shS$ denote the kernel. This means that the extension defined by $c \in H^1(Y,\shS)$
comes from the extension given by $0 \ra \shK \ra \shT_1 \ra \O_Y \ra 0$,
and $\tilde{c} \in H^1(Y, \shK)$.
$\shK$ is semistable, since its degree is $0$ and $\mu_{\rm max}(\shS) \leq 0$.
But then the image of $c$ is $0$ in every quotient sheaf of $\shS$ with negative maximal slope,
which contradicts the assumptions. Therefore we may assume that $\mu_{\rm max} (\shS') < 0$.

We want to apply Lemma \ref{slopecompare} to $\shS \subset \shS' =\shT$.
Assume that $\shS$ and $\shS'$ have a common subsheaf occuring in both Harder-Narasimhan filtrations.
Then they have the same maximal destabilizing subsheaf $\shF=\shS_1 =\shT_1$, which has negative degree.
If $c$ comes from $\tilde{c} \in H^1(Y,\shF)$, then $\shF \subset \shF' \subseteq \shS'$
and $\mu(\shF)= \deg(\shF)/\rk (\shF) < \deg(\shF)/(\rk(\shF)+1) = \mu(\shF')$, which contradicts
the maximality of $\shF$.
Hence the image of $c$ in $H^1(Y, \shS/\shF)$ is not zero and we can mod out
$\shF$ as before.

Therefore we may assume that $\shS$ and $\shS'$ do not have any common
subsheaf in their Harder-Narasimhan filtrations.
Then Lemma \ref{slopecompare} yields that $\sigma_i \leq \tau_{i+1}$,
and all these numbers are $\leq 0$ and moreover $\tau_i <0$.
Lemma \ref{numkrit} applied to $\alpha_i=- \sigma_i$ and $\beta_i=-\tau_i$
yields that $\sum_{i=1}^r \sigma_i^2 \geq \sum_{i=1}^{r+1} \tau_i^2$, and $>$ holds since $\tau_1 \neq 0$.
\qed

\begin{remark}
Suppose that $\shS$ is a semistable locally free sheaf of negative degree,
and let $c \in H^1(Y, \shS)$ with corresponding extension $\shS'$.
Then Theorem \ref{hilbertkunzaffinecriterion} together with Lemma \ref{semistableminimal} yield 
the inequalities
$$ \frac{\deg(\shS)^2}{r+1} \leq   \mumu (\shS') \leq  \frac{\deg(\shS)^2}{r} \, .$$
If $\shS'$ is also semistable, then we have equality on the left.
\end{remark}

\section{A Hilbert-Kunz criterion for solid closure}

We come now back to our original setting of interest,
that of a two-dimen\-sional normal standard-graded
domain $R$ over an algebraically closed field $K$.
A homogeneous $R_+$-primary ideal $I=(f_1 \komdots f_n)$ gives rise to the
syzygy bundle $\Syz(f_1 \komdots f_n)(0)$ on $Y= \Proj R$
defined by the presenting sequence
$$0 \lra \Syz(f_1 \komdots f_n)(m) \lra \bigoplus_{i=1}^n \O_Y(m-d_i)
\stackrel{f_1 \komdots f_n}{\lra} \O_Y(m) \lra 0 \, .$$
Another homogeneous element $f$ of degree $m$ yields
an extension
$$0 \lra \Syz(f_1 \komdots f_n)(m) \lra \Syz(f_1\komdots f_n,f)(m) \lra \O_Y \lra 0$$
which corresponds to the cohomology class $\delta(f) \in H^1(Y,\Syz(f_1 \komdots f_n)(m) )$
coming from the presenting sequence
via the connecting homomorphism
$$\delta: H^0(Y, \O_Y(m))=R_m \ra H^1(Y, \Syz(f_1 \komdots f_n)(m)) \, .$$
The Hilbert-Kunz multiplicities of the ideals and the Hilbert-Kunz slopes of the syzygy bundles
are related in the following way.

\begin{lemma}
\label{hilbertkunzmultiplicityslope}
Let $K$ denote an algebraically closed field of characteristic $0$.
Let $R$ denote a standard-graded two-dimensional normal $K$-domain, $Y= \Proj R$.
Let $I$ be a homogeneous $R_+$-primary ideal and let $f$ denote a homogeneous element of degree $m$.
Then the Hilbert-Kunz multiplicities $e_{HK}(I)= e_{HK}((I,f))$ are equal if and only
if the Hilbert-Kunz slopes of the corresponding syzygies bundles
$\mumu(\Syz(f_1 \komdots f_n)(m)) =\mumu(\Syz(f_1 \komdots f_n,f )(m))$ are equal.
\end{lemma}

\proof
Let $\mu_k$ and $r_k$ ($\tilde{\mu}_k$ and $\tilde{r}_k$)
denote the ranks and the slopes in the Harder-Narasimhan filtration
of $\Syz(f_1 \komdots f_n)(0)$ (of $\Syz(f_1 \komdots f_n, f)(0)$
respectively).
For the Hilbert-Kunz multiplicities of the ideals $(f_1 \komdots f_n)$ and $(f_1 \komdots f_n, f)$
we have to compare
$$e_{HK}(I) = \frac{1}{2 \deg(Y)} \big(\sum_{k=1}^t r_k \mu_k^2 - \deg(Y)^2 \sum_{i=1}^n d_i^2 \big)$$
and
$$e_{HK}((I,f))  =  \frac{1}{2\deg(Y)}
\big(\sum_{k=1}^{\tilde{t}} \tilde{r}_k \tilde{\mu}_k^2 - \deg(Y)^2 (m^2+\sum_{i=1}^n d_i^2) \big) \, .$$
The extension defined by $c=\delta(f)\in H^1(Y, \Syz(f_1 \komdots f_n)(m))$
is
$$0 \lra \shS=\Syz(f_1 \komdots f_n)(m) \lra \shS'=\Syz(f_1 \komdots f_n,f)(m) \lra \O_Y \lra 0$$
and the Hilbert-Kunz slopes of these sheaves are
due to Proposition \ref{hilbertkunzslopeproperties} (v)
(since $\deg(\Syz(f_1 \komdots f_n)(0)) = -\deg(Y) \sum_{i=1}^n d_i$)
$$\mumu(\shS)= \sum_{k=1}^t r_k \mu_k^2 +2 (- \sum_{i=1}^n d_i\deg(Y)) m \deg(Y) + (n-1)m^2 \deg(Y)^2$$
and $\mumu(\shS')=$
\begin{eqnarray*}
 &=& \sum_{k=1}^{\tilde{t}} \tilde{r}_k \tilde{\mu}_k^2
+ 2 (- ( \sum_{i=1}^n d_i +m) \deg(Y)) m \deg(Y) + nm^2 \deg(Y)^2  \cr
&=& \sum_{k=1}^{\tilde{t}} \tilde{r}_k \tilde{\mu}_k^2
- 2 (\sum_{i=1}^n d_i) m \deg(Y)^2 + (n-1)m^2 \deg(Y)^2 - m^2 \deg(Y)^2 \, .
\end{eqnarray*}
So  the difference is in both cases (up to the factor $1/2 \deg(Y)$)
$$ \sum_{k=1}^{\tilde{t}} \tilde{r}_k \tilde{\mu}_k^2 -   \sum_{k=1}^t r_k \mu_k^2 - m^2 \deg(Y)^2 \, .$$
Therefore $e_{HK}(I)= e_{HK}((I,f))$
if and only if
$$\mumu (\Syz(f_1 \komdots f_n)(m))= \mumu (\Syz(f_1 \komdots f_n,f)(m)) \, .$$
\qed

\begin{remark}
Let $0 \ra \shS \ra \shT \ra \shQ \ra 0$ denote a short exact sequence of locally free sheaves.
Then the alternating sum of the Hilbert-Kunz slopes,
that ist $\mumu(\shS) - \mumu(\shT) + \mumu(\shQ)$
does not changes when we tensor the sequence with an invertible sheaf.
This follows from Proposition \ref{hilbertkunzslopeproperties}(v).
For an extension $0 \ra \shS \ra \shS' \ra \O_Y \ra 0$
this number is $\geq 0$ by Theorem \ref{hilbertkunzaffinecriterion},
and we suspect that this is true in general.
From the presenting sequence
$0 \ra \Syz(f_1 \komdots f_n)(0) \ra \bigoplus_{i=1}^n \O(-d_i) \ra \O_Y \ra 0$
it follows via
$e_{HK}(I)= \frac{1}{2 \deg(Y)}( \mumu (\Syz(f_1 \komdots f_n)(0))- \mumu (\bigoplus_{i=1}^n \O(-d_i))$
that the Hilbert-Kunz multiplicity of an ideal is always nonnegative.
In fact $I=R$ is the only ideal with $e_{HK}(I)=0$. This follows
from Theorem \ref{hilbertkunzsolid} below, since $1 \not\in I^*$ for $I \neq R$.
\end{remark}

We come now to the main result of this paper.
Recall that the solid closure of an $\fom$-primary ideal $I=(f_1 \komdots f_n)$
in a two-dimensional normal excellent domain $R$
is given by the condition that
$f \in (f_1 \komdots f_n)^*$ if and only
$D(\fom) \subset \Spec R[T_1 \komdots T_n]/(f_1T_1 \plusdots f_nT_n+f)$ is not an affine scheme.
In positive characteristic this is the same as tight closure, see \cite[Theorem 8.6]{hochstersolid}.
In the case of an $R_+$-primary homogeneous ideal in a standard-graded normal $K$-domain this is equivalent
to the property that the torsor $ \PP({\shS '} ^\dual) - \PP(\shS^\dual)$ over the corresponding curve
$Y= \Proj R$ is not affine (see \cite[Proposition 3.9]{brennertightproj}). This relates solid closure to the setting of the previous section.

\begin{theorem}
\label{hilbertkunzsolid}
Let $K$ denote an algebraically closed field.
Let $R$ denote a standard-graded two-dimensional normal $K$-domain.
Let $I$ be a homogeneous $R_+$-primary ideal and let $f$ denote a homogeneous element.
Then $f \in I^*$ if and only if $e_{HK} (I)= e_{HK}((I,f))$.
\end{theorem}
\proof
If the characteristic is positive then this is a standard result from tight closure theory
as mentioned in the introduction.
So suppose that the characteristic is $0$.
Let $I=(f_1 \komdots f_n)$ be generated by homogeneous elements, and set $m =\deg(f)$.
The containment in the solid closure, $f \in (f_1 \komdots f_n)^*$, is equivalent with the non-affineness
of the torsor $ \PP({\shS '} ^\dual) - \PP(\shS^\dual)$ \cite[Proposition 3.9]{brennertightproj},
where $\shS = \Syz(f_1 \komdots f_n)(m)$
and $S'$ is the extension given by the cohomology class $\delta(f)$. Hence the result follows from
Theorem \ref{hilbertkunzaffinecriterion}
and Lemma \ref{hilbertkunzmultiplicityslope}.
\qed

\bibliographystyle{plain}

\bibliography{bibliothek}

\end{document}